\documentclass[11pt]{article}
\usepackage[centertags]{amsmath}
\usepackage{a4,amssymb,amsfonts,theorem,stmaryrd}
\usepackage[all]{xy}

\allowdisplaybreaks

\theorembodyfont{\slshape}
\newtheorem{thm}{Theorem}[section]
\newtheorem{prop}[thm]{Proposition}

\newtheorem{lemma}[thm]{Lemma}
\newtheorem{cor}[thm]{Corollary}
\theorembodyfont{\rmfamily} 
\newtheorem{Def}[thm]{Definition}
\newtheorem{rem}[thm]{Remark}

\newenvironment{pf}{\textit{Proof.}}{\hfill q.e.d.}

\newcommand{\R}{\mathbb{R}}
\newcommand{\N}{\mathbb{N}}
\newcommand{\Z}{\mathbb{Z}}
\newcommand{\C}{\mathbb{C}}

\newcommand{\cp}{cp}
\newcommand{\Ccinf}{C_{\cp}^{\infty}}

\newcommand{\E}{\mathcal{E}}

\newcommand{\cone}{\operatorname{cone}}
\newcommand{\tube}{\operatorname{tube}}

\newcommand{\M}{\mathbf{M}^3_{\kappa}}
\newcommand{\Hyp}{\mathbf{H}}
\newcommand{\Sph}{\mathbf{S}}
\newcommand{\Euc}{\mathbf{E}}

\newcommand{\dom}{\operatorname{dom}}

\newcommand{\im}{\operatorname{im}}

\renewcommand{\min}{min}
\renewcommand{\max}{max}

\newcommand{\Isom}{\operatorname{Isom^+}}
\newcommand{\isom}{\mathfrak{isom^+}}
\newcommand{\g}{\mathfrak{g}}

\newcommand{\SL}{\operatorname{SL}}

\newcommand{\dev}{\operatorname{dev}}
\newcommand{\hol}{\operatorname{hol}}

\newcommand{\tr}{\operatorname{tr}}

\newcommand{\res}{\operatorname{res}}

\newcommand{\supp}{\operatorname{supp}}

\newfont{\hugemath}{cmsy10 scaled 3000}

\begin{document}

\title{Global rigidity of 3-dimensional cone-manifolds}
\author{Hartmut Weiss}
\maketitle

\begin{abstract}
\noindent We prove global rigidity for compact hyperbolic and spherical
cone-3-manifolds with cone-angles $\leq \pi$ (which are not Seifert fibered in
the spherical case), furthermore for a class of hyperbolic cone-3-manifolds of finite volume with cone-angles $\leq \pi$, possibly with boundary consisting of totally geodesic hyperbolic turnovers. To that end we first generalize the local rigidity result
contained in \cite{Wei} to the setting of hyperbolic cone-3-manifolds of finite volume as above. We then use the
techniques developed in \cite{BLP} to deform the cone-manifold structure to a
complete non-singular or a geometric orbifold structure, where global rigidity holds due to
Mostow-Prasad rigidity, cf.~\cite{Mos}, \cite{Pra}, in the hyperbolic case, resp.~\cite{deR}, cf.~also \cite{Rot}, in the spherical case. This strategy has already been implemented
successfully by \cite{Koj} in the compact hyperbolic case if the singular
locus is a link using Hodgson-Kerckhoff local rigidity, cf.~\cite{HK}.
\end{abstract}


\section{Introduction}
Let $X$ be a compact, orientable hyperbolic (resp.~spherical) cone-3-manifold with cone-angles $\leq
\pi$. Let $\Sigma=\cup_{i=1}^N e_i$ be the singular locus and $M=X \setminus
\Sigma$ the smooth part of $X$, let further $(\alpha_1, \ldots, \alpha_N)$ be
the vector of cone-angles. Due to the above assumptions on the cone-angles, a
component of $\Sigma$ will either be a (connected) trivalent graph or a circle embedded geodesically into $X$. $M$ carries
a smooth Riemannian metric of sectional curvature $-1$ in the hyperbolic case, resp.~$+1$
in the spherical case, which is necessarily incomplete (the metric completion
of $M$ being given by glueing the singular locus back in).

Recall that cone-3-manifolds of curvature $\kappa \in \R$ with cone-angles $\leq2\pi$ are complete metric length spaces with curvature bounded from below by $\kappa$ in the triangle comparison sense. The local structure of a cone-3-manifold $X$ of curvature $\kappa$ is fixed by prescribing local models: For any point $x \in X$ the metric ball of sufficiently small radius centered at $x$ is required to be isometric to a truncated cone of curvature $\kappa$ over a spherical cone-surface, which will be called the (spherical) link of $x$. We will always require the links to be homeomorphic to the 2-sphere, which implies that the metric space $X$ is indeed homeomorphic to a 3-manifold. The singular locus $\Sigma$  is the union of those points in $X$, whose link is {\em not} isometric to the standard round 2-sphere $\Sph^2$. Similarly, a cone-surface of curvature $\kappa \in \R$ with cone-angles $\leq 2\pi$ is a metric space with curvature bounded from below by $\kappa$ which is locally isometric to a truncated cone of curvature $\kappa$ over a circle of length $\alpha$ with $\alpha \leq 2\pi$. The number $\alpha$ will be called the cone-angle.

By the classification of compact oriented spherical cone-surfaces with
cone-angles $\leq \pi$, the link of a singular point in a cone-3-manifold is
either a spherical turnover $\Sph^2(\alpha, \beta,\gamma)$ with $\alpha +
\beta + \gamma > 2\pi$ or a spherical football $\Sph^2(\alpha,\alpha)$. Here $\Sph^2(\alpha, \beta,\gamma)$ is the double of a spherical triangle with angles $\frac{\alpha}{2}, \frac{\beta}{2}, \frac{\gamma}{2}$ and $\Sph^2(\alpha, \alpha)$ is the double of a spherical bigon with (necessarily equal) angles $ \frac{\alpha}{2}, \frac{\alpha}{2}$. The spherical cone-manifold structures on $\Sph^2(\alpha,\beta,\gamma)$ and $\Sph^2(\alpha, \alpha)$ are unique once the cone-angles are fixed. This enumeration of possible links of singular points explains the particular structure of the singular locus if the cone-angles are $\leq \pi$.

For a basic introduction to the geometry of cone-3-manifolds and
cone-surfaces we refer the reader to \cite{CHK}, see also the introductory
sections of \cite{BLP}. The main motivation for studying cone-3-manifolds
comes from Thurston's approach to geometrization of
3-orbifolds: cone-3-manifolds provide a way to deform geometric orbifold
structures. The orbifold
theorem has now been proven in full generality by M.~Boileau, B.~Leeb and
J.~Porti, cf.~\cite{BLP}.

Following \cite{Por} we will say that a cone-3-manifold $X$ is {\em Seifert
  fibered} if $X$ carries a Seifert fibration such that the components of
  $\Sigma$ are leaves of the fibration. In particular $\Sigma$ is a union of
  circles and $M=X\setminus \Sigma$ is a Seifert fibered 3-manifold.

In \cite{Wei} we determine the local structure of the deformation space of hyperbolic (resp.~spherical) cone-manifold structures in the above setting. We find that
the space of cone-manifold structures is locally parametrized by the cone-angles, more precisely we have:
\begin{thm}[local rigidity]\label{local_rigidity_cp} 
Let $X$ be a compact, orientable hyperbolic (resp.~spherical) cone-3-manifold
with cone-angles $\leq\pi$ (which is not Seifert fibered in the spherical case). Then assigning the vector of cone-angles
$(\alpha_1,\ldots,\alpha_N)$, where $N$ is the number of edges contained in $\Sigma$, to a hyperbolic (resp.~spherical) cone-manifold
structure yields a local parametrization of the space of such structures near the given one.
\end{thm}
In particular we obtain:
\begin{cor}
Let $X$ be a hyperbolic (resp.~spherical) cone-3-manifold as in Theorem \ref{local_rigidity_cp}. Then there are no deformations of the hyperbolic (resp.~spherical) cone-manifold structure leaving the cone-angles fixed. 
\end{cor}
This complements an earlier theorem of
C.~Hodgson and S.~Kerckhoff, which states the same for compact, orientable hyperbolic cone-3-manifolds with singular locus a link, i.e.~a union of circles, but
where cone-angles $\leq 2\pi$ are allowed, cf.~\cite{HK}.

Following $\cite{Koj}$ we will say that two compact hyperbolic (resp.~spherical) cone-3-manifolds $X_1$
and $X_2$ have the same {\em topological type} if there exists a homeomorphism
of pairs $\phi:(X_1,\Sigma_1) \rightarrow (X_2,\Sigma_2)$. 

Local rigidity asserts that given a compact hyperbolic (resp.~spherical)
cone-3-manifold $X$ with cone-angles $\leq \pi$, then any compact hyperbolic
(resp.~spherical) cone-3-manifold $X'$ of the same topological type, which
after identifying the pairs $(X',\Sigma')$ and $(X,\Sigma)$ is {\em close}  to $X$ and
has the same cone-angles, will be isometric to $X$. Here being close to $X$ means
being close in the topology of the deformation space of $X$ as in \cite{HK} or \cite{Wei}.

We will say that a compact hyperbolic (resp.~spherical) cone-3-manifold $X$ is {\em globally
  rigid} if whenever $X'$ is another compact hyperbolic (resp.~spherical)
cone-3-manifold of the same topological type such that the cone-angles
around corresponding edges of $\Sigma$ and $\Sigma'$ coincide, then $X$ and $X'$ are isometric.

We will also be concerned with a class of hyperbolic cone-3-manifolds
of finite volume with cone-angles $\leq \pi$, namely those which have at most finitely many ends, all of
which are (smooth or singular) cusps with compact cross-sections, and which possibly
have totally geodesic turnover boundary. By the finiteness result of \cite{BLP}, cf.~also Corollary \ref{finiteness}, a hyperbolic
cone-3-manifold of finite volume, possibly with totally geodesic turnover boundary, satisfies these assumptions automatically if
the cone-angles are $\leq \eta < \pi$.

By the classification of compact oriented Euclidean cone-surfaces with
cone-angles $\leq \pi$, a singular cusp cross-section is either
a Euclidean turnover $\Euc^2(\alpha, \beta, \gamma)$ with $\alpha+\beta+\gamma
= 2\pi$ or a cone-surface of type $\Euc^2(\pi,\pi,\pi,\pi)$. Here $\Euc^2(\alpha, \beta, \gamma)$ is the
double of a Euclidean triangle with angles $\frac{\alpha}{2}, \frac{\beta}{2}, \frac{\gamma}{2}$.
Note that the Euclidean cone-manifold structure on $\Euc^2(\alpha, \beta,
\gamma)$ is unique (up to scaling). $\Euc^2(\pi,\pi,\pi,\pi)$ denotes a
Euclidean cone-manifold structure on $S^2$ with four cone-points of cone-angle
$\pi$. Note that such a structure is not uniquely determined, there is
actually a 2-dimensional family (up to scaling) of such structures. Smooth cusps are rank-2
cusps, i.e.~based on $T^2$. A totally geodesic boundary turnover will be a
hyperbolic turnover $\Hyp^2(\alpha,\beta,\gamma)$ with $\alpha+\beta+\gamma <
2\pi$, which is the double of a hyperbolic triangle with angles $\frac{\alpha}{2}, \frac{\beta}{2}, \frac{\gamma}{2}$. The hyperbolic
cone-manifold structure on $\Hyp^2(\alpha,\beta,\gamma)$ is unique.

For reasons which will become apparent later, we allow a little more
flexibility in the definition of topological type:
We say that two hyperbolic cone-3-manifolds of finite volume $X_1$
and $X_2$ as above have the same topological type, if after
truncating singular cusps (if there are any) and removing singular balls
around the vertices (if there are any), there exists a homeomorphism of pairs
$\phi:(\bar{X_1},\Sigma_1 \cap \bar{X_1}) \rightarrow (\bar{X_2},\Sigma_2
\cap \bar{X_2})$, where $\bar{X_1}$ and $\bar{X_2}$ denote the results of the
above operations on $X_1$ and $X_2$ respectively. Note however that if both
$X_1$ and $X_2$ are compact (without boundary), this coincides with the original definition.

For an extension of the local rigidity result to the setting of hyperbolic cone-3-manifolds of finite volume as above cf.~Theorem \ref{local_deformation}, Corollary \ref{local_rigidity} and Corollary \ref{hyperbolic_Dehn_surgery}, which require as an additional assumption that
$\Euc^2(\pi, \pi, \pi, \pi)$ doesn't occur as cusp cross-section.

We will say that a hyperbolic cone-3-manifold $X$ of finite volume as
above is globally rigid if whenever $X'$ is another hyperbolic cone-3-manifold of the same topological type such that the cone-angles
around corresponding edges of $\Sigma$ and $\Sigma'$ coincide, then $X$ and $X'$ are isometric.\\
\\
The main result in the hyperbolic case is:
\begin{thm}[global rigidity: compact case]\label{global_rigidity_hyperbolic}
Let $X$ be a compact, orientable hyperbolic cone-3-manifold with cone-angles $\leq \pi$. Then $X$ is globally rigid.
\end{thm}
This generalizes a theorem of S.~Kojima, which asserts
that compact, orientable hyperbolic cone-3-manifolds with cone-angles $\leq \pi$ and singular locus a link are globally rigid, cf.~\cite{Koj}.\\
\\
The same proof yields global rigidity for a class of hyperbolic cone-3-manifolds of finite volume:
\begin{thm}[global rigidity: finite-volume case]\label{global_rigidity_hyperbolic_finite_volume} 
Let $X$ be an orientable hyperbolic cone-3-manifold of finite volume with cone-angles $\leq\pi$, at most finitely many ends which are (smooth or singular) cusps with compact cross-sections and possibly boundary consisting of totally geodesic hyperbolic turnovers. Suppose that $\Euc^2(\pi,\pi,\pi,\pi)$ doesn't occur as cusp cross-section. Then $X$ is globally rigid.
\end{thm}
Note that by the finiteness result of \cite{BLP}, cf.~also Corollary \ref{finiteness}, a
hyperbolic cone-3-manifold of finite volume and cone-angles $\leq \eta < \pi$
has at most finitely many ends and all of them are cusps with compact (smooth
or singular) cross-sections. Therefore if the cone-angles are $\leq \eta <\pi$, Theorem \ref{global_rigidity_hyperbolic_finite_volume} deals with the general case of a finite-volume hyperbolic cone-3-manifold, possibly with totally geodesic turnover boundary.\\
\\
Apart from the extension of the local rigidity result of \cite {Wei}, cf.~also Theorem \ref{local_rigidity_cp}, to the setting of hyperbolic cone-3-manifolds of finite volume, cf.~Theorem \ref{local_deformation}, Corollary \ref{local_rigidity} and Corollary \ref{hyperbolic_Dehn_surgery}, we find worth stating the following corollaries of the proof of the main result:
\begin{cor}[smooth parts]\label{smooth_parts}
Let $X$ be an orientable hyperbolic cone-3-manifold of finite volume with cone-angles $\leq\pi$, at most finitely many ends which are (smooth or singular) cusps with compact cross-sections and possibly boundary consisting of totally geodesic hyperbolic turnovers. Suppose that $\Euc^2(\pi,\pi,\pi,\pi)$ doesn't occur as cusp cross-section. Then the smooth part of $X$ carries a complete (non-singular) hyperbolic structure of finite volume, possibly with totally geodesic boundary consisting of thrice-punctured spheres. 
\end{cor}
The following is a generalization of (a part of) the stability theorem of
\cite{BLP}, for the notion of {\em geometric convergence} used by \cite{BLP}
see also Definition \ref{geometric_convergence}:

\begin{cor}[topological stability] Let $(X_n,x_n)_{n \in \N}$ be a sequence of
  pointed orientable hyperbolic cone-3-manifolds of finite volume with cone-angles $\in [\zeta,\eta]$, $0 < \zeta \leq \eta < \pi$, possibly with boundary consisting of totally geodesic hyperbolic turnovers. Suppose that the $X_n$ have the same topological type and that there is uniform upper volume bound $vol(X_n) \leq V$. If the sequence $(X_n,x_n)_{n \in \N}$ converges geometrically to a pointed orientable hyperbolic cone-3-manifold $(X_{\infty}, x_{\infty})$, then $X_{\infty}$ has the same topological type as the $X_n$.
\end{cor}
\begin{pf}
Since $vol(X_{\infty}) \leq V$, the finiteness result of \cite{BLP}, cf.~also
Corollary \ref{finiteness}, implies that $X_{\infty}$ has at most finitely
many ends, all of which are cusps with compact (smooth or singular)
cross-sections. By Corollary \ref{smooth_parts}, the smooth part of
$X_{\infty}$ carries a complete (non-singular) hyperbolic structure of finite
volume, possibly with totally geodesic boundary consisting of thrice-punctured
spheres. Now the argument of the proof of Theorem 7.1 in \cite{BLP} applies.
\end{pf}\\
\\
In both corollaries we don't need to assume (as it is done in \cite{BLP}), that the underlying space is a small orbifold.\\
\\
The main result in the spherical case is:
\begin{thm}[global rigidity: spherical case]\label{global_rigidity_spherical} Let $X$ be a compact, orientable
  spherical cone-3-manifold with cone-angles $\leq \pi$. Suppose that $X$ is not Seifert
  fibered. Then $X$ is globally rigid.
\end{thm}
The strategy of the proof of the main result in the hyperbolic as in the spherical case is similar to the one applied in \cite{Koj}, nevertheless the presence of vertices in the singular locus causes some additional complications:

In the hyperbolic case we wish to construct a continuous, angle-decreasing family of hyperbolic cone-3-manifolds of constant topological type connecting the initially given cone-manifold structure with a complete (non-singular) hyperbolic structure of finite volume, possibly with totally geodesic boundary consisting of thrice punctured spheres. By Mostow-Prasad rigidity, cf.~\cite{Mos}, \cite{Pra}, applied to the double along the totally geodesic boundary, this structure is unique. We will assign cone-angle $0$ to the non-singular structure. We wish to conclude with a local rigidity theorem for the type of hyperbolic cone-3-manifolds through which we deform (resp.~with Thurston's hyperbolic Dehn surgery theorem applied to the double at cone-angle $0$) that the family of cone-manifold structures is unique if we prescribe cone-angles. 

Even if the initial hyperbolic cone-3-manifold is compact, we have to deform through a more general class of cone-3-manifolds, namely we have to allow a change in the geometry of the links: Recall that the sum of the cone-angles of three edges meeting at a vertex is always $>2\pi$. This is because the link of a vertex is a spherical turnover. If the cone-angles can be decreased such that this sum reaches $2 \pi$, the link necessarily becomes a horospherical Euclidean turnover in a singular cusp and the vertex disappears at infinity. Even further decreasing cone-angles amounts to deforming the link into a totally geodesic hyperbolic turnover, i.e.~the sum of the cone-angles of the three edges will be $<2\pi$. 

To carry out this strategy, two ingredients are essential: We need a local deformation theory for hyperbolic cone-3-manifolds with totally geodesic turnover boundary and singular cusps with Euclidean turnover cross-section. This generalizes the results contained in \cite{Wei}. It presents no further difficulty to also allow the presence of smooth rank-2 cusps here. Secondly we need to study degenerations of angle-decreasing families of hyperbolic cone-3-manifolds. Here we can apply the geometric results of \cite{BLP}, which underly their proof of the orbifold theorem.

In the spherical case we want to construct a continuous, angle-increasing family of spherical cone-3-manifolds of constant topological type connecting the initially given one with a compact spherical 3-orbifold with cone-angle $\pi$ along each edge. A spherical structure on a compact 3-orbifold is unique by a result of G.~de Rham, cf.~\cite{deR}, cf.~also \cite{Rot}. Since we are increasing cone-angles now, the issue of a change of geometry of the links doesn't arise here.

The author would like to thank Bernhard Leeb and Joan Porti for valuable
comments and suggestions. The results contained in this article rely strongly on the geometric results obtained by M.~Boileau, B.~Leeb and J.~Porti in \cite{BLP}. The author would further like to thank Steven Kerckhoff for useful discussions during the preparation of this article at Stanford University. This present work was in part funded by DFG-SPP 1154 "Globale Differentialgeometrie".

\section{Local rigidity of hyperbolic cone-3-manifolds of finite volume}

Let $X$ be a hyperbolic cone-3-manifold of finite volume with cone-angles
$\leq\pi$, at most finitely many ends which are (smooth or singular) cusps
with compact cross-sections and possibly boundary consisting of totally
geodesic hyperbolic turnovers. For the analytic results contained in this
section we may allow $\Euc^2(\pi,\pi,\pi,\pi)$ as cusp cross-section, only
later in the discussion of the local structure of the representation variety
do we have to exclude it.

In the above setting, the horospherical cross-section of each cusp is either a
flat torus $T^2=\R^2/\Gamma$, or a Euclidean turnover
$\Euc^2(\alpha,\beta,\gamma)$ with $\alpha,\beta,\gamma \leq \pi$ and $\alpha
+ \beta + \gamma = 2\pi$, resp.~a cone-surface of type
$\Euc^2(\pi,\pi,\pi,\pi)$. In the first case we are looking at a smooth rank-2 cusp, in the second we have three, resp.~four singular edges exiting the cusp.
Note that the Euclidean metric on the horospherical cross-section is only determined up to scaling. Recall that in a cusp neighbourhood the hyperbolic metric looks as follows:
$$
g^{hyp} = dt^2 + e^{-2t} g^{eucl}\,,
$$ 
where $g^{euc}$ is the Euclidean metric on the horospherical cross-section and $t \in[0,\infty)$.
A hyperbolic turnover, which appears as a component of the totally geodesic boundary, will be denoted by $\Hyp^2(\alpha,\beta,\gamma)$, where $\alpha, \beta, \gamma \leq \pi$ and $\alpha+\beta+\gamma<2\pi$. For the purposes of analysis we will attach the end
$$
[0,\infty) \times \Hyp^2(\alpha,\beta,\gamma) 
$$
with metric
$$
g^{hyp}=dt^2 + \cosh(t)^2 g^{\Hyp^2(\alpha,\beta,\gamma)}
$$
in the obvious way to the boundary, which is totally geodesic on both
sides. Let us denote the extended hyperbolic cone-3-manifold by $\hat{X}$, i.e.
$$
\hat{X} = X \cup_{\Hyp_i^2(\alpha,\beta,\gamma)} [0,\infty) \times \Hyp_i^2(\alpha,\beta,\gamma)\,.
$$ 
The extended hyperbolic cone-3-manifold $\hat{X}$ is now boundaryless but has infinite
volume ends instead. Clearly $\hat{X}$
  is homeomorphic to $X \setminus \partial X$.
Let $M$ denote the smooth part of $\hat{X}$ in the following.\\
\\
We consider the bundle of infinitesimal
isometries $\E = \mathfrak{so}(TM) \oplus TM$. The hyperbolic metric $g^{hyp}$
on $M$
yields a metric $h^{\E}$ on $\E$. The bundle $\E$ is furthermore equipped with
a flat connection which is given by
$$
\nabla_Y^{\E}(A,X)=(\nabla_YA-R(Y,X),\nabla_YX-AY)\,,
$$
where $\nabla$ denotes the Levi-Civita connection and $R$ the curvature tensor
of $g^{hyp}$. Parallel sections $\sigma=(A,X)$ correspond to Killing fields, more
precisely we have $\nabla^{\E}(A,X)=0$ if and only if $X$ is a Killing field and
$A=\nabla X$. $\E$ is in fact a bundle of Lie algebras, the Lie algebra
structure being given by
$$
[(A,X),(B,Y)]=([A,B]-R(X,Y),AY-BX)\,.
$$
Note that the metric $h^{\E}$ is {\em not} parallel with respect
to $g^{hyp}$, one rather has
$$
(\nabla^{\E}_X h^{\E})(\sigma,\tau) = -2h^{\E}(ad(X)\sigma,\tau)\,.
$$
Alternatively, the bundle $\E$ may be described as an associated bundle
$$
\E = \tilde{M} \times_{\pi_1 M} \isom \Hyp^3
$$
via the representation $Ad \circ \hol$, where $\hol: \pi_1M \rightarrow \Isom
\Hyp^3$ is the holonomy representation associated with the hyperbolic
structure on $M$. In particular, since $\isom \Hyp^3=\mathfrak{sl}_2(\C)$, the bundle $\E$ carries a parallel complex structure. For details we ask the reader to
consult \cite{Wei}.\\
\\
The flat connection $\nabla^{\E}$ extends via the Leibniz rule to an operator 
$$
d:\Omega^{\bullet}(M,\E) \rightarrow \Omega^{\bullet+1}(M,\E)\,,
$$
using its formal adjoint defined with respect to the metrics $g$ and $h^{\E}$
$$
d^t :  \Omega^{\bullet+1}(M,\E) \rightarrow \Omega^{\bullet}(M,\E)
$$
we may form the Hodge-Dirac operator $d+d^t$ and the Hodge-Laplace operator $\Delta = D^2 = dd^t + d^td$.\\
\\
If $\{e_1,e_2,e_3\}$ is a local orthonormal frame and $\{e^1,e^2,e^3\}$ the
dual coframe, we define
$$
\mathcal{D}:=\sum_{i=1}^3
\varepsilon(e^i)\nabla_{e_i}\;\text{and}\;T:=\sum_{i=1}^3\varepsilon(e^i)ad(e_i)\,,
$$
the formal adjoints are then given by
$$
\mathcal{D}^t=-\sum_{i=1}^3 \iota(e^i)\nabla_{e_i}\;\text{and}\;T^t=\sum_{i=1}^3\iota(e_i)ad(e_i)\,.
$$
Here $\varepsilon(e^i)$ denotes exterior multiplication with $e^i$ and $\iota(e_i)$ denotes interior multiplication with $e_i$.
Let further
$\Delta_\mathcal{D}:=\mathcal{D}\mathcal{D}^t+\mathcal{D}^t\mathcal{D}$ and
$H:=TT^t+T^tT$. Then we have the following Weitzenb\"ock formula, which is essentially due to
Y.~Matsushima and S.~Murakami, cf.~\cite{MM}, cf.~also \cite{HK}:
\begin{prop}
$\Delta = \Delta_{\mathcal{D}} + H$ and $H$ is strictly positive on 1-forms, i.e.~there
  exists a constant $c>0$ such that $(H\omega,\omega)_x\geq c(\omega,\omega)_x$ for all $\omega \in
  \Omega^1(M,\E)$ and $x\in M$.
\end{prop}
For $\omega \in \Omega_{cp}^1(M,\E)$, Stokes' theorem therefore yields 
$$
\int_M (\Delta \omega,\omega) = \int_M |\mathcal{D}|^2 + \int_M
(H\omega,\omega) \geq C \int_M |\omega|^2
$$
for some constant $C>0$, i.e.~$\Delta \geq C$ on compactly supported 1-forms.

\subsection{Analytic properties of cone-3-manifolds}

\noindent To define $L^2$-cohomology we consider the following subcomplex of
the de-Rham complex:
$$
\Omega^i_{L^2}(M,\E) = \{\omega \in \Omega^i(M,\E) : w \in L^2 \,\text{and}\,
dw \in L^2\}\,,
$$
which we will refer to as the smooth $L^2$-complex. Then
$H^{\bullet}_{L^2}(M,\E)$ is
by definition the cohomology of the smooth $L^2$-complex.

We generally think of a differential operator $P$ as acting on smooth,
compactly supported sections of some vector bundle. These constitute a dense
domain of definition inside the space of $L^2$-sections. We consider the
following closed extensions:
$$
\dom P_{\max} = \{ s \in L^2 : Ps \in L^2 \} 
$$
and 
$$
\dom P_{\min} = \{ s \in L^2: \exists s_n \in \Ccinf \text{ with } s_n
\overset{L^2}{\rightarrow} s,\, s_n \overset{L^2}{\rightarrow} Ps \} \,.
$$   
The maximal extension $P_{\max}$, resp.~the minimal extension $P_{\min}$ is then given by restricting $P$ (acting in a
distributional sense on $L^2$-sections) to $\dom P_{\max}$, resp.~$\dom P_{\min}$. By definition
$\dom P_{\max}$ is the maximal subspace of $L^2$ which is mapped to $L^2$
under $P$, whereas $\dom P_{\min}$ is just the closure of $\dom P =
\Ccinf$ with respect to the graph-norm of $P$. For details we ask the reader to consult $\cite{Wei}$.

We consider $d_{\max}$, the maximal
extension of $d$. It is not hard to see that $d^i_{\max}(\dom d^i_{\max}) \subset \dom
  d^{i+1}_{\max}$ and $d^{i+1}_{\max} \circ d^i_{\max} = 0$. Therefore the
$d_{max}$-complex
$$
\ldots \longrightarrow \dom d^i_{\max} \overset{d^i_{\max}}{\longrightarrow}
\dom d^{i+1}_{\max} \longrightarrow \ldots
$$
is defined as a complex in the sense of homological algebra. This is a
particular instance of a Hilbert-complex in the sense of \cite{BL}. Let us denote the cohomology of the $d_{\max}$-complex by $H_{\max}^{\bullet}$.
We define the $d_{\max}$-harmonic $i$-forms to be
$$
\mathcal{H}_{\max}^i = \ker d_{\max}^i\cap\ker(d^{i-1})^t_{\min}\,.
$$ 
We remind the reader of a regularity result:
\begin{thm}
The inclusion $\Omega_{L^2}^i(M,\E) \hookrightarrow \dom d^i_{\max}$ induces
an isomorphism in cohomology $H^i_{L^2}(M,\E) \cong  H^i_{\max}$.
\end{thm}
and a basic Hodge-type decomposition theorem:
\begin{thm}\label{weakHodge}
For each $i$
there is an orthogonal decomposition
$$
L^2(\Lambda^iT^*M\otimes\E)=\mathcal{H}_{\max}^i\oplus\overline{\im d^{i-1}_{\max}}\oplus\overline{\im (d^i)^t_{\min}}\,.
$$
\end{thm}
Both theorems in this form are due to J.~Cheeger, cf.~\cite{Ch1} and \cite{Ch2},
a more general discussion using the framework of Hilbert-complexes may be found in \cite{BL}. 

We define reduced $L^2$-cohomology
$$
h^i_{L^2}(M,\E):=\ker d^i_{max}/ \im \overline{\im d^{i-1}_{\max}}
$$
such that by the above theorem we have $\mathcal{H}^i_{max} \cong
h^i_{L^2}(M,\E)$. Furthermore the natural map $H^i_{max} \rightarrow
h^i_{L^2}(M,\E)$ is an isomorphism if and only if $d^{i-1}_{max}$ has closed range.

The principal aim of this section is to establish the following theorem, which extends the corresponding result in \cite{Wei} to the current setting:

\begin{thm}\label{L2vanishing} Let $(M,g^{hyp})$ be the smooth part of the
  extended hyperbolic cone-3-manifold $\hat{X}$ and $(\E,\nabla^{\E},h^{\E})$
  the flat bundle of infinitesimal isometries. If the cone-angles are $\leq \pi$, then $h^1_{L^2}(M,\E)=0$.
\end{thm}
In the compact case the range of $d^i_{max}$ is closed for all $i$ since $D_{min}=D_{max}$ has discrete spectrum, cf.~\cite{Wei}. Hence by the above remarks Theorem \ref{L2vanishing} implies $H^1_{L^2}(M,\E)=0$ in the compact case.

Note that essential selfadjointness of $D$ implies that $\Delta(d_{\max}) = \Delta_F$, cf.~\cite{Wei}. Here $\Delta_F$ denotes the Friedrichs
extension of $\Delta$. Therefore, taking into account the
Weitzenb\"ock formula and the fact that the Friedrichs extension preserves
lower bounds, the proof of Theorem \ref{L2vanishing} reduces to the
following analytic statement:
\begin{prop}\label{ess_self}
$D$ is essentially selfadjoint on $M$.
\end{prop}
\begin{pf}
Proving essential selfadjointness of $D$ amounts to proving that the minimal
and the maximal extension of $D$ coincide. Let $\omega \in \dom D_{\max}$ be
given. Without loss of generality we may assume $\omega$ to be smooth. We have
to show that already $\omega \in \dom D_{\min}$, i.e.~that there exists a
sequence of smooth, compactly supported forms $\omega_n$ satisfying
$$
\| \omega - \omega_n\|_{L^2} \underset{n \rightarrow \infty}{\longrightarrow} 0
$$
and
$$
\| D\omega - D\omega_n\|_{L^2} \underset{n \rightarrow \infty}{\longrightarrow} 0\,.
$$
Said in another way, the forms $\omega_n$ have to approximate $\omega$ with
respect to the graph-norm of $D$.\\
\\
We proceed in two steps:\\

\noindent Step 1: Approximation by forms with bounded support\\
\\
Here we essentially use the Gaffney cut-off argument in the
ends of $M$, cf.~\cite{Gaf}. Due to the comparatively simple geometry of the
ends, this can be made very explicit. We consider a cusp end
$$
E = [0,\infty) \times N
$$
with metric
$$
g^{hyp}=dt^2 + e^{-t^2} g^{eucl}\,,
$$
where $(N,g^{eucl})$ denotes the (smooth part of the) horospherical cross-section,
which may be $T^2$, $\Euc(\alpha,\beta,\gamma)$ or $\Euc(\pi,\pi,\pi,\pi)$,
or an infinite volume end
$$
E = [0,\infty) \times \Hyp^2(\alpha,\beta,\gamma) 
$$
with metric
$$
g^{hyp}=dt^2 + \cosh(t)^2 g^{\Hyp^2(\alpha,\beta,\gamma)}
$$
in the following. We consider cut-off functions
$\varphi_n=\varphi_n(t)$, which will be chosen suitably later, and we set
$\omega_n=\varphi_n\omega$. Recall that the Hodge-Dirac operator satisfies the
following product rule:
$$
D(\varphi_n\omega) = c(d\varphi_n)\omega + \varphi_n D\omega\,,
$$
where $c(d\varphi_n)=\varepsilon(d\varphi_n) - \iota (d\varphi_n)$ denotes Clifford
multiplication with $d\varphi_n$. In particular we obtain
$$
D\omega - D \omega_n = (1-\varphi_n)D\omega - c(d\varphi_n)\omega \,.
$$
If one has a sequence of cut-off functions
$\varphi_n$ such that $\supp \varphi_n$ exhaust $M$, then
$$
\| \omega - \omega_n\|_{L^2} = \| (1-\varphi_n) \omega \|_{L^2} \underset{n \rightarrow \infty}{\longrightarrow} 0
$$
and
$$
\| (1-\varphi_n)D\omega\|_{L^2} \underset{n \rightarrow \infty}{\longrightarrow} 0\,.
$$
It remains to arrange that $\sup_M | d\varphi_n| \rightarrow 0$ as $n \rightarrow
\infty$, then also
$$
\|c(d\varphi_n)\omega\|_{L^2} \leq \sup_M |d\varphi_n| \cdot \|\omega\|_{L^2} \underset{n \rightarrow \infty}{\longrightarrow} 0\,.
$$ 
We can construct such a sequence of cut-off functions as follows:

Let $n \geq 3$. Define $\varphi_n|_{[0,n]} \equiv 1$,
$\varphi_n|_{[2n,\infty)} \equiv 0$. On $[n+1,2n-1]$ let $\varphi_n$ linearly
  decay from $1-1/n$ to $1/n$ with slope $-1/n$. Finally interpolate smoothly
  on $[n,n+1]$ and $[2n-1,2n]$ with $|d\varphi_n| \leq 1/n$. The functions
  $\varphi_n$ extend by $1$ to $M$. Clearly $\supp \varphi_n$ exhaust $M$ and
  $|d\varphi_n|\leq 1/n$.\\
 
\noindent Step 2: Approximation by forms with compact support\\
\\
In view of Step 1 it is sufficient to approximate a form with bounded support
by forms with compact support. But here the arguments given in \cite{Wei} for
the hyperbolic case - using the separation of variables techniques developed
by J.~Br\"uning and R.~Seeley in \cite{BS} - directly apply.\\
\\
This finishes the proof of Proposition \ref{ess_self} and of Theorem
\ref{L2vanishing}. 
\end{pf}
\begin{prop}\label{imageL2cohomology}
$\im (H^1_{L^2}(M,\E) \rightarrow H^1(M,\E))=\{0\}$.
\end{prop}
\begin{pf}
Let $\omega \in \Omega^1_{L^2}(M,\E)$ with $d \omega =0$. We need to show that there exists $s \in \Gamma(M,\E)$ with $ds=\omega$. 

We decompose $\omega$ according to Theorem \ref{weakHodge}, i.e.
$$
L^2\Lambda^1 T^*M \otimes \E = \mathcal{H}^1_{max} \oplus \overline{\im d^0_{max}} \oplus \overline{\im (d^1)^t_{min}}\,,
$$
and use Theorem \ref{L2vanishing} and the fact that $\overline{\im (d^1)^t_{min}}=(\ker d^1_{max})^\perp$ to conclude that $\omega \in \overline{\im d^0_{max}}$. Hence there exists a sequence $s_n \in \dom d^0_{max}$ such that
$$
\|d s_n - \omega\|_{L^2}  \underset{n \rightarrow \infty}{\longrightarrow} 0\,,
$$
where we can assume without loss of generality that the $s_n$ are smooth.

Let now $\Omega \subset M$ be a relatively compact domain with $\partial \Omega$ smooth such that $\bar{\Omega}$ is a compact core for $M$. In particular we have
$$
H^0(\Omega, \E) = \{ s \in \Gamma(\Omega,\E) : \nabla s =0 \} =\{0\}
$$
since the holonomy representation is irreducible, cf.~Section \ref{repvar}. This implies that there exists a constant $C>0$ such that 
$$
\|\sigma\|_{L^2(\Omega)} \leq C \cdot \| \nabla \sigma\|_{L^2(\Omega)}
$$
for all $\sigma \in W^{1,2}(\Omega)$. Namely, if we had a sequence $\sigma_n \in W^{1,2}(\Omega)$ with $\lim_{n \rightarrow \infty}\|\nabla\sigma_n\|_{L^2(\Omega)}=0$ and $\|\sigma_n\|_{L^2(\Omega)}=1$, then by compactness of the embedding $W^{1,2}(\Omega) \hookrightarrow L^2(\Omega)$ after passing to a subsequence $\sigma_n$ would converge in $L^2(\Omega)$ and weakly in $W^{1,2}(\Omega)$ to a section $\sigma \in W^{1,2}(\Omega)$. Clearly $\|\sigma\|_{L^2(\Omega)}=1$ and $\nabla \sigma =0$, which contradicts $H^0(\Omega,\E)=\{0\}$.

Going back to our original sequence $s_n$, we may now use the above estimate together with Rellich's lemma to pass to a subsequence such that $s_n\vert_{\Omega}$ converges in $L^2(\Omega)$ to some section $s_{\Omega} \in W^{1,2}(\Omega)$. Finally, by using an exhaustion of $M$ by domains $\Omega_m$ as above, we obtain a diagonal subsequence $s_n$ which satisfies
$$
s_n \overset{L^2_{loc}}{\longrightarrow} s \in L^2_{loc} \quad \text{ and } \quad \nabla s_n \overset{L^2}{\longrightarrow} \omega \,.
$$
It remains to show that  $s$ is smooth and that $ds=\omega$. For $\varphi \in C_{cp}^{\infty}$ we find some $m \in \N$ such that $\supp \varphi \in \Omega_m$. Now
$$
\int_M (s,d^t\varphi) = \lim _{n \rightarrow \infty} \int_{\Omega_m} (s_n,d^t \varphi) = \lim_{n \rightarrow \infty} \int_{\Omega_m} (ds_n, \varphi) = \int_M(\omega,\varphi)\,,
$$
which shows that $ds=\omega$ in a distributional sense. Smoothness of $s$ is a consequence of elliptic regularity since $\Delta s =d^t \omega$ is smooth.   
\end{pf}\\
\\
An immediate corollary, cf.~\cite{HK}, cf.~also \cite{Wei}, is the following:
\begin{cor}\label{halfdimensional}
Let $\bar{M}$ be a compact core of the smooth part of the extended hyperbolic cone-3-manifold $\hat{X}$. Then:
\begin{enumerate}
\item The natural map $H^1(\bar{M},\E) \rightarrow H^1(\partial \bar{M},\E)$ is
injective.
\item $\dim H^1(\bar{M},\E)=\frac{1}{2}\dim H^1(\partial \bar{M},\E)$.
\end{enumerate}
\end{cor}

\subsection{The representation variety}\label{repvar}
Let $X$ be a hyperbolic cone-3-manifold of finite volume with cone-angles $\leq\pi$, at most finitely many ends which are (smooth or singular) cusps with compact cross-sections and possibly boundary consisting of totally geodesic hyperbolic turnovers. We will assume from now on that $\Euc^2(\pi,\pi,\pi,\pi)$ does not occur as cusp cross-section.

A compact core $\bar{M}$ may be obtained from $M$ by the
following procedure: Fix $\varepsilon >0$ and $A>0$. We truncate the cusps at
Euclidean turnovers and tori of area $A$ and the infinite volume
ends at the totally geodesic hyperbolic turnovers. From what remains we remove
the $\varepsilon$-tube around $\Sigma$. We obtain a decomposition 
$$M = \bar{M} \cup U^{A,tg}_{\varepsilon}(\Sigma) \cup U^A\,,$$ 
where $U^A$ is the union of the nonsingular cusp neighbourhoods,
$U^{A,tg}_{\varepsilon}(\Sigma)$ is the union of the $\varepsilon$-tube around
$\Sigma$, the cusp neighbourhoods of the singular cusps and the infinite
volume ends beyond the totally geodesic turnovers, and $\bar{M}$ is the
compact core of $M$.\\
\\
We recall some basic facts about representation spaces of discrete groups:
Let $\Gamma$ be a finitely generated discrete group and let $G$ be a Lie
group. We denote by $R(\Gamma, G)$ the space of group homomorphisms $\rho:
\Gamma \rightarrow G$ equipped with the compact-open topology. By fixing a presentation $\langle \gamma_1, \ldots, \gamma_n \vert (r_i)_{i \in I} \rangle$ the {\em representation
variety} $R(\Gamma,G)$ may be identified with a subset of $G^n$. If $G$ is a
complex algebraic group, e.g.~if $G=\SL_2(\C)$, $R(\Gamma,G)$ acquires the
structure of a complex algebraic set. Let $X(\Gamma,G)$ be the quotient of
$R(\Gamma,G)$ under the natural action of $G$ by conjugation. We will equip
the space $X(\Gamma,G)$ with the quotient topology.

For
$\rho \in R(\Gamma,G)$ let further $H^i(\Gamma, \g)$ denote the $i$-th group
cohomology group with coefficients in the representation $Ad \circ \rho$. If the relations cut out $R(\Gamma,G)$ transversely, $T_{\rho}R(\Gamma,G)$ may be identified with the space of 1-cocycles $Z^1(\Gamma, \mathfrak{g})$. The tangent space to the orbit through $\rho$ is given by the space of 1-coboundaries $B^1(\Gamma,G)$. Hence, if $G$ acts properly and freely on $R(\Gamma,G)$, then $X(\Gamma,G)$ will be smooth and, if $\chi$ denotes the conjugacy class of a representation $\rho$, $T_{\chi}X(\Gamma, G)$ may be identified with $H^1(\Gamma,\mathfrak{g})$.

If $\Gamma = \pi_1M$, we can form the associated bundle
$\E = \tilde{M} \times_{\pi_1 M} \mathfrak{g}$. Since $M$ is a 3-manifold, the period map
\begin{align*}
P: H^1(M,\E) & \longrightarrow H^1(\pi_1M,\mathfrak{g})\\
[\omega] & \longmapsto [\gamma\mapsto\textstyle\int_{\gamma}\omega]
\end{align*}
is an isomorphism. 

The holonomy representation of a hyperbolic cone-manifold
structure
$$
\hol : \pi_1\bar{M} \rightarrow \Isom \Hyp^3
$$
may always be lifted to $\widetilde{\Isom \Hyp^3} = \SL_2(\C)$ by a theorem of
M.~Culler, cf.~\cite{Cul}. We will denote such a lift again by $\hol$.\\
\\
In the following let $\rho$ be the restriction of $\hol$ to $P$, where $P$ is either the smooth part of the spherical link of a vertex, the smooth part of a horospherical cross-section of a singular cusp (a {\em horospherical link}) or the smooth part of a totally geodesic hyperbolic turnover (a {\em hyperspherical link}). Topologically $P$ is a thrice punctured sphere, whose fundamental group is the free group on 2 generators. We use the following presentation:
$$
\pi_1 P = \langle \mu_1, \mu_2, \mu_3 \,\vert\, \mu_1\mu_2\mu_3 = 1 \rangle \,.
$$
We obtain
$$
R(\pi_1P,\SL_2(\C)) = \{ (A_1,A_2,A_3) \in \SL_2(\C)^3 \,\vert\, A_1A_2A_3 = 1 \}\,. 
$$
This realizes $R(\pi_1P,\SL_2(\C))$ as a complex codimension-3 submanifold in $\SL_2(\C)^3$, in particular $\dim_{\C}R(\pi_1P,\SL_2(\C))=6$.

\begin{lemma}\label{geometry_of_links}
Let $A_1, A_2, A_3 \in\Isom \Hyp^3$ elliptic satisfying $A_1A_2A_3=1$. If the $A_i$ do not fix the same line in $\Hyp^3$, then one of the following (mutually exclusive) cases holds:
\begin{enumerate}
\item The $A_i$ fix a unique point $x \in \Hyp^3$ and preserve the foliation of $\Hyp^3$ by distance spheres centered at $x$ leafwise.
\item The $A_i$ fix a unique point $x_{\infty} \in \partial_{\infty} \Hyp^3$ and preserve the foliation of $\Hyp^3$ by horospheres centered at $x_{\infty}$ leafwise.
\item The $A_i$ preserve a unique totally geodesic $\Hyp^2 \subset \Hyp^3$.
\end{enumerate} 
\end{lemma}
\begin{pf}
Since the $A_i$ are elliptic, they each fix a line $\gamma_i\in \Hyp^3$. Now either $\gamma_1=\gamma_2=\gamma_3$, or $\gamma_1, \gamma_2, \gamma_3$ are pairwise distinct. In the latter case, the $\gamma_i$ intersect in a unique point $x \in \Hyp^3$, $x_{\infty} \in \partial_{\infty}\Hyp^3$ or they don't intersect at all, in which case it is not hard to show that they are perpendicular to a unique totally geodesic $\Hyp^2 \subset \Hyp^3$. 
\end{pf}\\
\\
Let now $\rho$ be the restriction of the holonomy to the smooth part of a spherical link $\Sph^2(\alpha,\beta,\gamma)$. As in \cite{Wei} we obtain that $\rho$ is irreducible and therefore the centralizer of $\rho$ is trivial: $Z(\rho(\pi_1 P))=\{\pm 1\}$. It follows in particular that $H^0(\pi_1 P, \mathfrak{sl}_2(\C))=0$.
  
If $\rho$ is the restriction of the holonomy to the smooth part of a
horospherical link $\Euc^2(\alpha,\beta,\gamma)$, then $\rho$ is reducible, i.e.~the image of $\rho$ fixes a point in the boundary of $\Hyp^3$. Nevertheless, since the image of $\rho$ does not fix a line in $\Hyp^3$, we still have $Z(\rho(\pi_1 P))=\{\pm 1\}$ and $H^0(\pi_1 P, \mathfrak{sl}_2(\C))=0$.

Finally, if $\rho$ is the restriction of the holonomy to the smooth part of a hyperspherical link $\Hyp^2(\alpha,\beta,\gamma)$, then again $\rho$ will be irreducible, hence $Z(\rho(\pi_1 P))=\{\pm 1\}$ and $H^0(\pi_1 P,
\mathfrak{sl}_2(\C))=0$ also in this case.

We obtain that the holonomy representation $\hol: \pi_1 \bar{M}
\rightarrow \SL_2(\C)$ is irreducible. Namely if $X$ contains a spherical or
a hyperspherical link, then by the above said $\hol$ will be irreducible,
otherwise the argument in the proof of Lemma 6.35 in \cite{Wei}
applies. Furthermore this implies $Z(\hol(\pi_1\bar{M}))= \{\pm 1\}$ and $H^0(\bar{M},\E)=0$.

For $\rho \in  R(\pi_1P,\SL_2(\C))$ let $t_{\mu_i}(\rho)=\tr \rho(\mu_i)$ for $i=1,2,3$. Then we have:

\begin{lemma}
The differentials $\{dt_{\mu_1},dt_{\mu_2},dt_{\mu_3}\}$ are $\C$-linearly independent in $T^*_{\rho} R(\pi_1P,\SL_2(\C))$.
\end{lemma}
\begin{pf}
Taking into account that $Z(\rho(\pi_1 P))=\{\pm 1\}$ in all three cases, the proof given in \cite{Wei} for the case of spherical links works similarly for the case of horospherical and hyperspherical links.
\end{pf}\\
\\
As a consequence of the above lemma, together with $Z(\rho(\pi_1
P))=\{\pm 1\}$, we obtain using the implicit function theorem, that
$X(\pi_1P,\SL_2(\C))$ is smooth near $[\rho]$ and that
$\{t_{\mu_1},t_{\mu_2},t_{\mu_3}\}$ are local holomorphic coordinates, in
particular this implies $\dim_{\C} X(\pi_1P, \SL_2(\C))=3$ near $[\rho]$.

Let $N$ be the number of edges contained in $\Sigma$ and $\tau$ the number of the toral cusps. As in \cite{Wei} we obtain using that $H^0(\pi_1P,\mathfrak{sl}_2(\C))=0$ if $P$ is the smooth part of a (spherical, horospherical or hyperspherical) link together with the usual arguments for the toral cusps, cf.~\cite{Kap}:
\begin{lemma}
$H^1(\partial \bar{M}, \mathcal{E}) \cong \C^{2N+2\tau}$.
\end{lemma}
and hence using Corollary \ref{halfdimensional}:
\begin{cor}
$H^1(\bar{M}, \mathcal{E}) \cong \C^{N+\tau}$.
\end{cor}
Away from the vertices and the ends, the singular
tube $U^{A,tg}_{\varepsilon}(\Sigma)$ can be given Fermi-type coordinates $(r_i,\theta_i,z_i)$ with $r_i\in(0,\varepsilon)$, $\theta_i \in
\R/\alpha_i\Z$ and $z_i\in (0,l_i)$. Here $\alpha_i$ is the cone-angle around the
$i$-th edge and $l_i$ the length of a finite subsegment of the $i$-th edge, around which we define coordinates. Note
that edges exiting the ends have infinite length. The
hyperbolic metric is then given by:
$$
g^{hyp}=dr^2+\sinh^2(r)d\theta_i+\cosh^2(r)dz_i\,.
$$
We choose a function $\varphi_i=\varphi_i(z_i)$ such that $\varphi_i(0)=0$,
$\varphi_i(l_i)=l_i$ and
$d\varphi_i\vert_{(0,\delta)}=d\varphi_i\vert_{(l_i-\delta,l_i)}=0$ for
$\delta>0$. Then $d\varphi_i \in \Omega^1(U^{A,tg}_{\varepsilon}(\Sigma))$ is
well-defined and so are
\begin{align*}
\omega_{tws}^i &= d\varphi_i \varotimes \sigma_{\partial/\partial\theta_i}\\
\omega_{len}^i &= d\varphi_i \varotimes \sigma_{\partial/\partial z_i}\,.
\end{align*}  
Note that these forms are supported away from the vertices of the
singularity.
As in \cite{Wei} we easily obtain:
\begin{lemma}
The forms $\omega_{tws}^i$ and $\omega_{len}^i$ are bounded on
$U^{A,tg}_{\varepsilon}(\Sigma)$ and hence in particular $L^2$. 
\end{lemma}
and using again that $H^0(\pi_1P,\mathfrak{sl}_2(\C))=0$ if $P$ is the smooth part of a (spherical, horospherical or hyperspherical) link:
\begin{lemma}
The de-Rham cohomology classes of the differential forms 
$$\{ \omega_{tws}^1, \omega_{len}^1, \ldots , \omega_{tws}^N,
\omega_{len}^N \}$$ 
are $\R$-linearly independent in $H^1(\partial \bar{M}, \E)$.
\end{lemma}
Note that with respect to the parallel complex structure on $\E$, one has $
\omega_{tws}^i=\sqrt{-1}\,\omega_{len}^i$, therefore the classes of the forms 
$\{\omega_{len}^1, \ldots , \omega_{len}^N \}$ are $\C$-linearly independent in $H^1(\partial \bar{M}, \E)$.

In the following let $\rho$ be the restriction of $\hol$ to $\partial \bar{M}$. Let $N$ be the number of edges contained in $\Sigma$.
Let $\tau$ be the number of the toral cusps and let us fix meridians, i.e.~simple closed curves, $m_1,\ldots,m_{\tau}$ on the corresponding boundary tori in $\partial \bar{M}$.

For $\rho \in  R(\pi_1 \partial \bar{M},\SL_2(\C))$ let $t_{\mu_i}(\rho)=\tr \rho(\mu_i)$ for $i=1,\ldots,N$, resp.~let $t_{m_i}(\rho)=\tr \rho(m_i)$ for $i=1,\ldots,\tau$. Then the glueing procedure described in \cite{Wei} together with the usual arguments for the toral cusps, cf.~\cite{Kap}, yields the following:
\begin{lemma}
$\rho$ is a smooth point in $R(\pi_1\partial \bar{M}, \SL_2(\C))$, furthermore the differentials $\{dt_{\mu_1}, \ldots ,dt_{\mu_N}\} \cup \{dt_{m_1},\ldots,dt_{m_{\tau}}\}$ are $\C$-linearly independent in $T^*_{\rho} R(\pi_1\partial \bar{M}, \SL_2(\C))$.
\end{lemma}
Here, since $\partial \bar{M}$ may be disconnected, $\pi_1\partial\bar{M}$ refers to the fundamental group of the one-point union of the connected components of $\partial \bar{M}$, i.e.~the free product of the fundamental groups of the components. If $\partial_1\bar{M}, \ldots, \partial_l\bar{M}$ is the collection of the connected components of $\partial\bar{M}$, we therefore have:
$$
R(\pi_1\partial \bar{M}, \SL_2(\C)) = \prod_{i=1}^l R(\pi_1 \partial_i\bar{M}, \SL_2(\C))\,.
$$
Let $\rho_i$ be the restriction of $\rho$ to the $i$-th boundary component. As in \cite{Wei}, resp.~\cite{Kap} for the toral cusps, one shows that $X(\pi_1 \partial_i \bar{M},\SL_2(\C))$ is smooth near $\chi_i=[\rho_i]$. The functions $t_{\mu_i}$, resp.~$t_{m_i}$, descend to functions
$$
t_{\mu_i}: \prod_{i=1}^l X(\pi_1 \partial_i\bar{M}, \SL_2(\C)) \rightarrow \C
$$
resp.
$$
t_{m_i}: \prod_{i=1}^l X(\pi_1 \partial_i\bar{M}, \SL_2(\C)) \rightarrow \C\,.
$$
The differentials $\{dt_{\mu_1}, \ldots ,dt_{\mu_N}\} \cup \{dt_{m_1},\ldots,dt_{m_{\tau}}\}$ remain $\C$-linearly independent.
We have a natural restriction map:
$$
\res: X(\pi_1\bar{M},\SL_2(\C)) \rightarrow \prod_{i=1}^l X(\pi_1 \partial_i\bar{M}, \SL_2(\C))\,.
$$
Lemma \ref{halfdimensional} implies that $\res$ is locally an immersion around
$\chi$, furthermore $X(\pi_1\bar{M},\SL_2(\C))$ has $\C$-dimension $N+\tau$, $\prod_{i=1}^l X(\pi_1 \partial_i\bar{M}, \SL_2(\C))$ has $\C$-dimension $2N+2\tau$.

Without loss of generality we may assume that $\partial_1\bar{M}, \ldots, \partial_k\bar{M}$ are those components of $\partial \bar{M}$, which bound components of the singular tube. Consequently $\partial_{k+1}\bar{M}, \ldots, \partial_{k+\tau}\bar{M}$ will be those tori in $\partial \bar{M}$, which bound rank-2 cusps. We consider the level sets
$$
V=\{t_{\mu_1}\equiv t_{\mu_1}(\chi), \ldots,
t_{\mu_N}\equiv t_{\mu_N}(\chi) \} \subset \prod_{i=1}^k X(\pi_1 \partial_i\bar{M}, \SL_2(\C))
$$
and
$$
W=\{t_{m_1}\equiv t_{m_1}(\chi), \ldots,
t_{m_{\tau}}\equiv t_{m_{\tau}}(\chi) \}  \subset \prod_{i=k+1}^{k+\tau} X(\pi_1 \partial_i\bar{M}, \SL_2(\C))\,.
$$
Note that locally around $\chi$ the level sets $V$ and $W$ are
half-dimensional submanifolds, furthermore $\prod_{i=1}^k X(\pi_1 \partial_i\bar{M}, \SL_2(\C))$ has $\C$-dimension $2N$ and $\prod_{i=k+1}^{k+\tau} X(\pi_1 \partial_i\bar{M}, \SL_2(\C))$ has $\C$-dimension $2\tau$. As in \cite{Wei} we obtain:
\begin{lemma}
The cohomology classes
of the cocycles 
$\{z_{len}^1,\ldots,z_{len}^N\}$ provide a $\C$-basis of $T_{\chi}V$.
\end{lemma}
Using the $L^2$-vanishing theorem as in \cite{Wei} together with the usual arguments for the toral cusps, cf.~\cite{Kap}, we obtain that  $V \times W$ and the image of $\res$ in $\prod_{i=1}^{l} X(\pi_1 \partial_i\bar{M}, \SL_2(\C))$ meet transversally. We have therefore proved the following:
\begin{thm}\label{main}
Let $\bar{M}$ be a compact core of the extended hyperbolic cone-3-manifold $\hat{X}$. Let $\{\mu_1, \ldots, \mu_N\} \cup \{m_1,\ldots,m_{\tau}\}$ be the family of meridians. Then the map 
\begin{align*}
X(\pi_1 \bar{M},\SL_2(\C)) & \longrightarrow \C^{N+\tau}\\
\chi & \longmapsto (t_{\mu_1}(\chi),\ldots,t_{\mu_N}(\chi),t_{m_1}(\chi),\ldots,t_{m_{\tau}}(\chi) )
\end{align*}
is locally biholomorphic near $\chi=[\hol]$.
\end{thm}
These coordinates can be converted into a family of more geometric ones:

If we fix $m_i,l_i$ generators of $\pi_1T^2_i$, where $T^2_i$ is the horospherical
cross-section of the $i$-th toral cusp,
we can pass to Dehn surgery coefficients, cf.~\cite{CHK} or \cite{BP},
i.e.~we obtain a local parametrization
\begin{align*}
X(\pi_1 \bar{M},\SL_2(\C)) & \longrightarrow \C^{N}\times (\R^2 \cup
\{\infty\}/\pm 1)^{\tau}\\
\chi & \longmapsto
(t_{\mu_1}(\chi),\ldots,t_{\mu_N}(\chi),[x_1,y_1],\ldots,[x_{\tau},y_{\tau}] ) \,.
\end{align*}
For the $i$-th toral cusp, $[x_i,y_i]=\infty$ corresponds to the cusped structure,
whereas lines through the origin with rational slope correspond to (families of) cone-manifold structures on a fixed topological filling. 

More precisely, if
$x_i/y_i \in \mathbb{Q} \cup \{ \infty \}$, write $x_i/y_i=p_i/q_i$ with $p_i,
q_i$ coprime integers. The cone-manifold structure corresponding to
$[x_i,y_i]$ is obtained by glueing in a singular solid torus with cone-angle
$2\pi|p_i/x_i|$ such that the curve $p_im_i+q_il_i$ bounds the singular disk. 
Note in particular that moving out to $\infty$ along the line with slope
$p_i/q_i$ corresponds to monotonically decreasing the cone-angle towards $0$
while preserving the topological type of the cone-3-manifold.

Dehn surgery coefficients different from those above do not give rise to cone-manifold structures. We will generally say that a hyperbolic structure which is obtained by varying the $i$-th Dehn surgery coefficient has a {\em Dehn surgery type singularity} at the $i$-th toral cusp. 

If we restrict to characters corresponding to holonomies of cone-manifold structures, i.e.~meridians are mapped to pure rotations with rotation-angle given by the cone-angle, then we can recover the trace of meridian from the cone-angle, cf.~for example \cite{HK}, \cite{Wei}.

Since mapping a hyperbolic structure to the character corresponding to its
holonomy representation induces a local homeomorphism between the deformation
space of hyperbolic structures on $M$ and the space $X(\pi_1\bar{M},\SL_2(\C))$, cf.~\cite{Gol}, we obtain the following local deformation theorem:

\begin{thm}\label{local_deformation}
Let $X$ be an orientable hyperbolic cone-3-manifold of finite volume with cone-angles $\leq\pi$, at most finitely many ends which are (smooth or singular) cusps with compact cross-sections and possibly boundary consisting of totally geodesic hyperbolic turnovers. Suppose that $\Euc^2(\pi,\pi,\pi,\pi)$ doesn't occur as cusp cross-section. Then assigning the vector of cone-angles
$(\alpha_1,\ldots,\alpha_N)$ and the family of Dehn surgery coefficients $([x_1,y_1], \ldots, [x_{\tau},y_{\tau}])$ to a hyperbolic cone-manifold structure with Dehn surgery type singularities yields a local parametrization of the space of such structures near the given one. 
\end{thm}
The geometry of the link of a vertex is determined by the cone-angles $\alpha,
\beta$ and $\gamma$ of the adjacent edges, cf.~Lemma \ref{geometry_of_links}:

If $\alpha+\beta+\gamma>2\pi$, the link is spherical, i.e.~given by a spherical turnover $\Sph^2(\alpha, \beta,\gamma)$. The vertex is just an ordinary vertex sitting inside the cone-3-manifold.
If $\alpha + \beta + \gamma =2\pi$, the link is horospherical, i.e.~given by a horospherical Euclidean turnover $\Euc^2(\alpha,\beta,\gamma)$. The vertex itself should now be considered to sit at infinity. Finally, if $\alpha+\beta+\gamma<2\pi$, the link is hyperspherical, i.e.~given by a totally geodesic hyperbolic turnover $\Hyp^2(\alpha, \beta, \gamma)$. The vertex should now be considered to sit beyond the totally geodesic boundary.

By leaving the Dehn surgery coefficients untouched, we obtain a local rigidity result for a class of hyperbolic cone-3-manifolds of finite volume analogous to the one contained in \cite{Wei}, cf.~also Theorem \ref{local_rigidity_cp}:

\begin{cor}[local rigidity: finite-volume case]\label{local_rigidity}Let $X$ be an orientable hyperbolic cone-3-manifold of finite volume with cone-angles $\leq\pi$, at most finitely many ends which are (smooth or singular) cusps with compact cross-sections and possibly boundary consisting of totally geodesic hyperbolic turnovers. Suppose that $\Euc^2(\pi,\pi,\pi,\pi)$ doesn't occur as cusp cross-section. Then assigning the vector of cone-angles
$(\alpha_1,\ldots,\alpha_N)$  to a hyperbolic cone-manifold
structure yields a local para-metrization of the space of such structures near the given one.
\end{cor}
Similarly, by leaving the cone-angles untouched, we obtain the conclusion of Thurston's hyperbolic Dehn surgery theorem in the setting of hyperbolic cone-3-manifolds of finite volume:

\begin{cor}[hyperbolic Dehn surgery]\label{hyperbolic_Dehn_surgery} 
Let $X$ be an orientable hyperbolic cone-3-manifold of finite volume with cone-angles $\leq\pi$, at most finitely many ends which are (smooth or singular) cusps with compact cross-sections and possibly boundary consisting of totally geodesic hyperbolic turnovers. Suppose that $\Euc^2(\pi,\pi,\pi,\pi)$ doesn't occur as cusp cross-section. Then for each toral cusp the conclusion of Thurston's hyperbolic Dehn surgery theorem holds, in particular all but finitely many Dehn fillings (per cusp) are hyperbolic.
\end{cor}
Note that by Theorem \ref{local_deformation} we can moreover {\em independently} perturb the cone-angles and perform the operation of hyperbolic Dehn surgery at a toral cusp. This will be of importance in the proof of the main theorem.

\section{Global Rigidity}

\subsection{Geometric properties of cone-3-manifolds}

In this section we mostly review material contained in \cite{BLP}. We will state
results only as far as directly needed for our applications. For proofs and
more details we ask the reader to consult \cite{BLP}.

In the following let $X$ be a
cone-3-manifold of curvature $\kappa$. For $r>0$ (and $r \leq
\pi/\sqrt{\kappa}$ if $\kappa>0$) and $S$ a spherical cone-surface, 
the {\it standard ball} of radius $r$ with {\it link} $S$ is simply the
$\kappa$-cone over $S$, which we as usual denote by $\cone_{\kappa,(0,r)}S$.

For $p \in X$ consider $B_r(p)$, the embedded $r$-ball around p in $X$. Recall
that for
$\varepsilon >0$ small enough we have by definition of cone-manifold structure 
$$
B_{\varepsilon}(p) \cong \cone_{\kappa,(0,\varepsilon)}S(p) \,,
$$
where $S(p)$ is the link of $p$ in $X$.

The {\it injectivity radius} of $X$ at $p$, $r_{inj}(p)$, is the supremum over
all $r>0$ such that 
$$
B_{r}(p) \cong \cone_{\kappa,(0,r)}S(p) \,.
$$
The {\it cone-injectivity radius} of $X$ at $p$, $r_{cone-inj}(p)$, is the
supremum over all $r>0$ such that $B_r(p)$ is contained in an embedded
standard ball, i.e.~such that there exists $q \in X$ and $R>0$ with $B_r(p)
\subset B_R(q)$ and
$$
B_R(q) \cong \cone_{\kappa,(0,R)}S(q) \,.
$$
For $\varrho>0$, $X$ is said to be {\it $\varrho$-thick} (at a point $p$) if it contains
an embedded smooth standard ball of radius $\varrho$ (centered at $p$).

Let $X$ and $Y$ be metric spaces and $\varepsilon >0$. Following \cite{BLP} we call a map $\phi: X
\rightarrow Y$ a {\it $(1+\varepsilon)$-bilipschitz embedding} if
$$
(1+\varepsilon)^{-1}\cdot d(x_1,x_2) < d(\phi(x_1),\phi(x_2)) < (1+\varepsilon)\cdot d(x_1,x_2)
$$ 
holds for all $x_1, x_2 \in X$.
\begin{Def}[geometric convergence]{\rm\cite{BLP}}
\label{geometric_convergence}
Let $(X_n,x_n)_{n \in \N}$ be a sequence of pointed cone-3-manifolds. We say that the sequence $(X_n,x_n)$
converges {\it geometrically} to a pointed cone-3-manifold $(X_{\infty},x_{\infty})$
if for every $R>0$ and $\varepsilon>0$ there exists $N=N(R,\varepsilon) \in
\mathbb{N}$ such that for all $n \geq N$ there is a
$(1+\varepsilon)$-bilipschitz embedding $\phi_n: B_R(x_{\infty}) \rightarrow
X_n$ satisfying:
\begin{enumerate}
\item $d(\phi_n(x_{\infty}),x_n) < \varepsilon$,
\item $B_{(1-\varepsilon)R}(x_n) \subset \phi_n(B_R(x_{\infty}))$, and
\item $\phi_n(B_R(x_{\infty}) \cap \Sigma_{\infty}) = \phi_n(B_R(x_{\infty}))
  \cap \Sigma_n$.
\end{enumerate}
\end{Def}
Note that if the $X_n$ have curvature $\kappa_n \in \R$, then $X_{\infty}$
will have curvature
$\kappa_{\infty}=\lim_{n \rightarrow \infty}\kappa_n$. 
The cone-angle along an edge of $\Sigma_{\infty}$ will be the limit of the
cone-angles along the corresponding edge in $\Sigma_n$. Note however that part
of the singular locus of the approximating cone-3-manifolds may disappear at
the limit by going to infinity.

Given a geometrically convergent sequence $(X_n,x_n)$ we may without loss of
generality assume that the $(1+\varepsilon)$-bilipschitz embeddings $\phi_n$ restrict to smooth
$(1+\varepsilon)$-bilipschitz diffeomorphisms
$$
F_n : B_R(x_{\infty}) \cap M_{\infty} \rightarrow
\phi_n(B_R(x_{\infty})) \cap M_n \,,
$$
where $M_{\infty}$, resp.~$M_{n}$
denotes the smooth part of $X_{\infty}$, resp.~$X_n$. Let us assume in the
following that $\kappa_n=\kappa$ independent of $n$ and therefore also
$\kappa_{\infty}=\kappa$. We may further assume that (possibly after 
composing with isometries) 
$$
\dev_n \circ \widetilde{F}_n: \widetilde{B_R(x_{\infty}) \cap M_{\infty}} \rightarrow \M
$$ 
converges to $\dev_{\infty}$ restricted to $\widetilde{B_R(x_{\infty}) \cap
  M_{\infty}}$ uniformly on compact sets (together with all derivatives), and that
$$
\hol_n \circ (F_n)_*: \pi_1(B_R(x_{\infty}) \cap M_{\infty}) \rightarrow \Isom \M
$$
converges to $\hol_{\infty}$ restricted to $\pi_1(B_R(x_{\infty}) \cap
M_{\infty})$ in the compact-open topology.
\begin{thm}[compactness]{\rm \cite{BLP}}\label{compactness}
Let $(X_n,x_n)_{n \in \N}$ be a sequence of pointed cone-3-manifolds with curvatures $\kappa_n
\in [-1,1]$, cone-angles $\leq \pi$ and possibly with totally geodesic
boundary. Suppose that for some $\rho >0$
each $X_n$ is $\rho$-thick at $x_n$. Then (possibly after
passing to a subsequence) the sequence $(X_n,x_n)$ converges geometrically to a pointed cone-3-manifold
$(X_{\infty},x_{\infty})$ with curvature $\kappa_{\infty}=\lim_{n\rightarrow\infty} \kappa_n$.
\end{thm}
In the hyperbolic case, M.~Boileau, B.~Leeb and J.~Porti derive a thick-thin
decomposition in the spirit of the classical Margulis lemma applied to
complete hyperbolic 3-manifolds. To state the results we need some further
notation, which is introduced in \cite{BLP}:

Recall that an embedded connected surface $S$ in $\Hyp^3$ is called umbilic, if all principal
curvatures have the same value $pc(S)$. Its intrinsic curvature is then given
by $\kappa_S=-1 + pc(S)^2$. An umbilic surface $S$ is called spherical if $\kappa_S >0$,
horospherical if $\kappa_S=0$ and hyperspherical if $\kappa_S<0$.
 
Recall further that $\Hyp^3$ carries the following natural foliations by umbilic
surfaces:
\begin{enumerate}
\item The foliation by distance spheres of a point $p \in \Hyp^3$: the leaves
  are spherical. 
\item The foliation by horospheres centered at $p \in
  \partial_{\infty}\Hyp^3$: the leaves are horospherical.
\item The foliation by surfaces equidistant to a totally geodesic $\Hyp^2
  \subset \Hyp^3$: the leaves are hyperspherical.
\end{enumerate}
Let $S$ be a cone-surface
of curvature $\kappa_S \in \{-1,0,1\}$. Depending on the curvature $\kappa_S$
we define the {\it complete tube over $S$}, $\tube_{-1}S$, as follows:
\begin{enumerate}
\item If $\kappa_S=1$, then $\tube_{-1}S$ is just the complete hyperbolic cone over $S$:
$$
\tube_{-1} S = \cone_{-1}S \,.
$$ 
\item If $\kappa_S=0$, then $\tube_{-1}S$ is the complete hyperbolic cusp with horospherical
cross-section $S$:
$$
\tube_{-1}S = \mathbb{R} \times S
$$
with metric
$$
g^{hyp} = dt^2 + e^{-2t}g^S \,.
$$
\item If $\kappa_S=-1$, then $\tube_{-1}S$ is the complete hyperbolic neck with totally geodesic
central leaf $S$:
$$
\tube_{-1}S = \mathbb{R} \times S
$$
with metric
$$
g^{hyp} = dt^2 + \cosh(t)^2g^S \,.
$$
\end{enumerate}
In any case, $\tube_{-1}S$ carries a natural foliation by umbilic surfaces
equidistant to $S$. Its leaves are
spherical if $\kappa_{S}=1$, horospherical if $\kappa_{S}=0$ and
hyperspherical if $\kappa_{S}=-1$.

Following \cite{BLP}, an {\it umbilic tube} will be a closed connected subset of
$\tube_{-1}S$, which is a union of leaves of the natural umbilic foliation.
\begin{thm}[thin part]{\rm \cite{BLP}}\label{thin_part}
For $D>0$ and $0<\zeta\leq\eta<\pi$ there exist constants $i=i(D,\zeta,\eta)>0$, $P=P(\zeta,D)>0$, $\rho=\rho(D,\zeta,\eta)>0$ such that the following holds:
If $X$ is an orientable hyperbolic cone-3-manifold (without boundary) with cone-angles $\in[\zeta,\eta]$ and $diam(X) \geq D$, then $X$ contains a (possibly
empty) disjoint union $X^{thin}$ of submanifolds which belong to the following
list:
\begin{enumerate}
\item smooth Margulis tubes, i.e.~tubular neighbourhoods of closed geo-desics
  and smooth cusps of rank one or two,
\item tubular neighbourhoods of closed singular geodesics,
\item umbilic tubes with turnover cross-sections which have strictly convex boundary, i.e.~standard (singular) balls, cusps and necks.
\end{enumerate} 
Furthermore, the boundary of each component of $X^{thin}$ is non-empty,
strictly convex with principal curvatures $\leq \pi$ and each of its (at most
two) components contains a smooth point $p$ with $r_{inj}(p) \geq \rho$
(measured in $X$); each component of $X^{thin}$ contains an embedded smooth
standard ball of radius $\rho$; all singular vertices are contained in
$X^{thin}$, and on $X\setminus X^{thin}$ holds $r_{cone-inj} \geq i$.
\end{thm}
Not that the thin part is canonically foliated: smooth cusps by the
horospherical cross-sections, tubular neighbourhoods of (smooth or singular)
geodesics by the distance tori and umbilic tubes by the natural umbilic foliations. 

\begin{cor}[thickness]{\rm \cite{BLP}}\label{thickness}
There exists $r=r(D,\zeta,\eta) > 0$ such that if $X$ is as in Theorem
\ref{thin_part}, then $X$ is $r$-thick, i.e.~$X$ contains an embedded smooth
standard ball of radius $r$.
\end{cor}

\begin{cor}[finiteness]{\rm \cite{BLP}}\label{finiteness}
Let $X$ be as in Theorem \ref{thin_part} and suppose in addition that
$vol(X)<\infty$. Then $X$ has finitely many ends and all of them are (smooth
or singular) cusps with compact cross-sections.
\end{cor}
By doubling along the boundary we obtain corresponding statements for
hyperbolic cone-3-manifolds with totally geodesic boundary.

\subsection{The hyperbolic case}
Let $X$ be a hyperbolic cone-3-manifold of finite volume with cone-angles $\leq\pi$, at most finitely many ends which are (smooth or singular) cusps with compact cross-sections and possibly boundary consisting of totally geodesic hyperbolic turnovers. Suppose that $\Euc^2(\pi,\pi,\pi,\pi)$ does not occur as cusp cross-section. We will describe how to obtain a continuous family of hyperbolic cone-3-manifolds of the same topological type connecting the given one with a complete (non-singular) hyperbolic 3-manifold of finite volume, possibly with boundary consisting of thrice-punctured spheres. 

Let $N$ denote the
number of edges contained in $\Sigma$ and $(\alpha_1, \ldots ,\alpha_N)$ the
vector of cone-angles. Using local rigidity we may assume that the cone-angles
$\alpha_i$ are already
strictly smaller than $\pi$, i.e.~$\alpha_i \leq \eta < \pi$.\\
\\ 
The argument will proceed in two steps:
\begin{enumerate}
\item Given a vector $(\beta_1, \ldots, \beta_N)$ of small target-angles,
  construct a continuous angle-decreasing family of hyperbolic
  cone-3-manifolds $(X_{\tau})_{\tau \in [0,1]}$ of constant topological type
  connecting $X=X_1$ with a hyperbolic cone-3-manifold $X_0$ with cone-angles $(\beta_1, \ldots, \beta_N)$.

\item Given a hyperbolic cone-3-manifold $X$ with small cone-angles, construct
  a non-singular hyperbolic structure on the smooth part of $X$ and connect it
  with $X$ by a continuous angle-decreasing family of hyperbolic cone-3-manifolds of constant topological type.
\end{enumerate}
\begin{rem}
A hyperbolic cone-3-manifold will be considered to have {\it small}
cone-angles if the cone-angles are small enough to force all
links to be hyperspherical, e.g.~cone-angles $<\frac{2\pi}{3}$ would be
sufficient to guarantee this.
\end{rem}

\subsubsection{Deforming into small cone-angles}
Let $N$ denote the number of edges contained in $\Sigma$ and let $(\alpha_1,
\ldots, \alpha_N)$ be the vector of cone-angles of $X$. Let $(\beta_1, \ldots,
\beta_N)$ with $0 < \beta_i \leq \alpha_i$ be a specified vector of target-angles. We
consider the interval
\begin{align*}
\mathcal{I}= \{ & t \in [0,1] : \exists \text{ a family of
  hyperbolic cone-3-manifolds } (X_{\tau})_{\tau \in [t,1]} \\
& \text{of constant topological type with } X_1=X \text{ and cone-angles }\\ 
& (\tau\alpha_1+(1-\tau)\beta_1, \ldots,\tau\alpha_N+(1-\tau)\beta_N) \}\,.
\end{align*}
$\mathcal{I}$ is nonempty since clearly $1 \in
\mathcal{I}$. We claim that $\mathcal{I}$ is open and closed, which implies that
$\mathcal{I}=[0,1]$, i.e.~there exists a continuous angle-decreasing family of hyperbolic cone-3-manifolds $(X_{\tau})_{\tau \in [0,1]}$ with constant topological type connecting $X=X_1$ with a hyperbolic cone-3-manifold $X_0$ with cone-angles $(\beta_1, \ldots, \beta_N)$.\\
\\
{\bf(openness)} Let $t \in \mathcal{I}$ and a corresponding family of
hyperbolic cone-3-manifolds $(X_{\tau})_{\tau \in [t,1]}$ be given. If $t>0$, we can further continuously decrease the cone-angles of $X_t$ using local rigidity.\hfill q.e.d.\\
\\
{\bf(closedness)} Let $(t_n)_{n \in \mathbb{N}}$ be a sequence in $\mathcal{I}$, w.l.o.g.~we may assume that $t_n$ is decreasing with $\lim_{n \rightarrow \infty}=t_{\infty}\in[0,1]$. We claim that $t_{\infty} \in \mathcal{I}$.
Corresponding to $t_n\in\mathcal{I}$ there exists a family of hyperbolic cone-3-manifolds $(X^n_{\tau})_{\tau \in [t_n,1]}$ with the above properties. Note that for $m>n$ the families $X^n_{\tau}$ and $X^m_{\tau}$ have to coincide on the interval $[t_n,1]$ due to local rigidity, such that we are effectively given a single family $(X_{\tau})_{\tau \in (t_{\infty},1]}$ defined now on a half-open interval.

Recall Schl\"afli's differential formula for the volumes of a smooth family of compact cone-3-manifolds of curvature $\kappa\in\mathbb{R}$, cf.~for example \cite{CHK}:
$$
\kappa \frac{d}{d\tau} vol(X_{\tau}) = \frac{1}{2} \sum_{i=1}^N l_i \frac{d \alpha_i}{d \tau}
$$
Here as usual $\alpha_i$ denotes the cone-angle around the $i$-th edge and $l_i$ the length of the $i$-th edge. If we have a smooth family of hyperbolic cone-3-manifolds of finite volume with at most finitely many {\em smooth} cusps of rank 2, such that in particular all edges have finite length, Schl\"afli's formula continues to hold without modification.
  
\begin{lemma}\label{lowerdiameterbound}
There exists a constant $D>0$ such that $diam(X_{\tau}) \geq D$ for $\tau \in (t_{\infty},t_n]$ and $n\in \mathbb{N}$ sufficiently large.
\end{lemma}
\begin{pf}
We choose $n$ large enough such that no horospherical links occur on the interval $(t_{\infty},t_n]$. Schl\"afli's formula applied to the doubles of $X_{\tau}$ along totally geodesic turnover boundaries implies that $vol(X_{\tau})$ is a (strictly) decreasing function of $\tau$. This
implies that $vol(X_{\tau}) \geq vol(X_{t_n})>0$ for $\tau \in (t_{\infty}, t_n]$. The lower curvature bound now implies a lower
diameter bound $diam(X_{\tau}) \geq D$ for some $D > 0$.
\end{pf}
\begin{lemma}\label{uppervolumebound}
There exists a constant $V < \infty$ such that $vol(X_{\tau}) \leq V$ for $\tau \in (t_{\infty},t_n]$ and $n \in \mathbb{N}$ sufficiently large.
\end{lemma}
\begin{pf}
Again we choose $n$ large enough such that no horospherical links occur on the
interval $(t_{\infty},t_n]$. Kojima's straightening argument,
  cf.~the proof of Proposition 1.3.2 in \cite{Koj}, easily generalizes to the setting of hyperbolic cone-3-manifolds of finite volume with at most finitely many {\em smooth} cusps of rank 2. Applied to the doubles of $X_{\tau}$ along totally geodesic turnover boundaries it yields the result.
\end{pf}\\
\\
We can now apply Lemma
  \ref{lowerdiameterbound} and Corollary \ref{thickness} to find a uniformly
  thick base point $x_n$ in $X_n:=X_{t_n}$ if $n \in \mathbb{N}$ is sufficiently large. Theorem \ref{compactness} provides
  us (after taking a subsequence, if necessary) with a geometric limit
  $(X_{\infty},x_{\infty})=\lim_{n\rightarrow \infty}(X_n,x_n)$. $X_{\infty}$ will be a hyperbolic cone-3-manifold, using Lemma
  \ref{uppervolumebound} we obtain that $vol(X_{\infty}) \leq V < \infty$.

We may now use the finiteness result of \cite{BLP} and the description of the
thin part to analyze the geometry of the limit and how it is attained.

First of all, the finiteness result of $\cite{BLP}$, cf.~also Corollary
\ref{finiteness}, says that $X_\infty$ has only finitely many ends, all of
which are cusps with compact (smooth or singular) horospherical
cross-sections. Since the cone-angles are $\leq \eta < \pi$, as singular
cross-section only $\Euc^2(\alpha, \beta,\gamma)$ can occur. In particular
this implies that $X_{\infty}$ has a compact core with horospherical boundary
consisting of tori and Euclidean turnovers.

Fix $R>0$ large enough such that $B_R(x_{\infty}) \subset X_{\infty}$ contains
a compact core $\bar{X}_{\infty}$
with horospherical boundary such that
$X_{\infty} \setminus \bar{X}_{\infty}$ is contained in the thin part of
$X_{\infty}$. Let $\bar{M}_{\infty}$ denote the
smooth part of $\bar{X}_{\infty}$. For $n \in \mathbb{N}$ large, the bilipschitz embeddings $\phi_n$ may be
slightly modified, such that for a component $H \subset \partial
\bar{X}_{\infty}$ the image $H_n:=\phi_n(H)$ is a leaf in the canonical
foliation of the thin part of $X_n$, cf.~Section 7.1 in \cite{BLP}.

The description of the thin parts of hyperbolic cone-3-manifolds contained in
\cite{BLP}, cf.~also Theorem
\ref{thin_part}, can be used to determine the geometry of the regions $X_n \setminus \phi_n(\bar{X}_{\infty})$.
Similar to Lemma 7.4 in \cite{BLP} we obtain in our situation:
\begin{lemma}\label{degenerations}
For $n \in \mathbb{N}$ sufficiently large, each component of $X_n \setminus \phi_n(\bar{X}_{\infty})$ is contained in the
thin part of $X_n$ and is either a singular ball, a (smooth or singular) solid
torus or a smooth rank-$2$ cusp.
\end{lemma}
\begin{pf}
First of all, we may choose $n \in \mathbb{N}$ large enough, such that $X_n$
has no horospherical links. Let $H \subset \partial
\bar{X}_{\infty}$ be a component, i.e.~a horospherical
cross-section of a (smooth or singular) cusp. If $H$ is a turnover, then $H_n:=\phi_n(H) \subset X_n$ will be
an umbilic turnover. Since $X_n$ has no horospherical links, $H_n$ cannot be
horospherical. Since the cone-angles are decreasing with $n$, the sum of
the cone-angles of a hyperspherical turnover in $X_n$ is bounded away from $2\pi$,
i.e.~for $n\in\mathbb{N}$ large, $H_n$ cannot be
hyperspherical. Therefore $H_n$ must be spherical and bounds a singular
ball according to Theorem \ref{thin_part}. 

If $H$ is a torus, then $H_n$ will either be a distance torus of a (smooth or
singular) geodesic and therefore
bounds a (smooth or singular) solid torus, or a horospherical cross-section of a smooth cusp and
therefore bounds a smooth rank-$2$ cusp according to Theorem \ref{thin_part}.
\end{pf}\\
\\
In other words, the only degenerations that can occur are tubes around short (smooth
or singular) closed geodesics opening into rank-2 cusps. Note in particular
that $X_{\infty}$ has potentially more rank-2 cusps than $X_n$ and potentially
less singular edges. 

Ultimately we want to
be able to rule out such degenerations as well. The argument used in
\cite{BLP} in the corresponding situation uses the fact, that the underlying
space of the $X_n$ is a small orbifold in their case. Since we don't want to
assume this, their argument doesn't apply here.

We obtain from Lemma \ref{degenerations} that $M_n$ arises topologically from $\bar{M}_{\infty}$ as a Dehn filling of some of the boundary tori in $\partial
\bar{M}_{\infty}$. We claim that this is also
geometrically true if $n$ is large enough. Namely as a consequence of geometric convergence we obtain
diffeomorphisms
$$
F_n : \bar{M}_{\infty} \rightarrow F_n(\bar{M}_{\infty}) \subset M_n
$$
such that (possibly after composing with isometries)
$$
\dev_n \circ \widetilde{F}_n: \widetilde{\bar{M}_{\infty}} \rightarrow \Hyp^3
$$ 
converges to $\dev_{\infty}$ restricted to $\widetilde{\bar{M}_{\infty}}$
uniformly on compact sets (together with all derivatives),
and that
$$
\hol_n \circ (F_n)_*: \pi_1(\bar{M}_{\infty}) \rightarrow \SL_2(\C)
$$
converges to $\hol_{\infty}$ restricted to $\pi_1(\bar{M}_{\infty})$ in the
compact-open topology.

This means that $F_n(\bar{M}_{\infty})$ may be viewed as a small deformation
of $\bar{M}_{\infty}$ via $F_n$ if $n$ is large enough, in particular it will be controlled by local
rigidity combined with hyperbolic Dehn surgery, cf.~Theorem \ref{local_deformation} and
Corollaries \ref{local_rigidity} and \ref{hyperbolic_Dehn_surgery}. 

Fixing a sufficiently large $n$ to make the above true, we can now use Theorem
\ref{local_deformation} applied to $\bar{M}_{\infty}$ to continuously decrease the cone-angles of
$F_n(\bar{M}_{\infty})$ to the limit angles 
$$(t_{\infty}\alpha_1+(1-t_{\infty})\beta_1,\ldots,t_{\infty}\alpha_N+(1-t_{\infty})\beta_N)$$
while keeping the Dehn surgery coefficients corresponding
to tori in $\partial \bar{M}_{\infty}$ mapping to distance tori of smooth
closed geodesics or to horospherical tori of smooth cusps constant. More
precisely, if a torus in $\partial \bar{M}_{\infty}$ is mapped to a distance
torus of a singular geodesic, let us say the $i$-th edge, then we can move
with the Dehn surgery coefficient corresponding to that torus on the line of
constant slope out to infinity and therefore decrease the cone-angle $\alpha_i$
as much as we like, in particular we can reach the limit angle
$t_{\infty}\alpha_i+(1-t_{\infty})\beta_i$. We can certainly leave the Dehn
surgery coefficients corresponding to tori mapping to distance tori of smooth
closed geodesics or to horospherical tori of smooth cusps constant, and finally
we can use the cone-angle parameters of Theorem \ref{local_deformation}
applied to $\bar{M}_{\infty}$ to decrease the cone-angles of the remaining
edges to the limit angles.

We have therefore extended the initially given family of hyperbolic cone-3-manifolds
$(X_{\tau})_{\tau \in (t_{\infty},1]}$ to the closed interval $[t_{\infty},1]$
and hence established closedness of $\mathcal{I}$.\hfill q.e.d.\\ 
\\
A posteriori we also find, using
again local rigidity, that the cusp-opening degenerations have not occurred at all.

\subsubsection{Deforming into cone-angle $0$}

Let us now assume that we are given a hyperbolic cone-3-manifold $X$ with
small cone-angles, i.e.~the links of all vertices are totally geodesic
hyperbolic turnovers. As usual let $N$ denote the number of edges contained in
$\Sigma$ and $(\alpha_1, \ldots ,\alpha_N)$ the vector of cone-angles. Let us
further assume for a moment that the singular locus does contain vertices, otherwise we are
anyway in a situation covered by Kojima's result, cf.~\cite{Koj}. 

We can now form the double of
$X$ along the totally geodesic boundary, let us denote it by
$2X$. By construction the singular locus of $2X$ has no vertices.
Therefore by a well-known procedure using Thurston's hyperbolization
theorem for Haken manifolds, cf.~Theorem 1.2.1 in \cite{Koj}, the smooth part of $2X$ carries a
complete (non-singular) hyperbolic structure of finite volume. Let us denote the smooth part of
$2X$ equipped with this structure by $2M_0$. We will see that
$2M_0$ is the double of a complete hyperbolic 3-manifold of finite volume
$M_0$ with totally geodesic
boundary consisting of thrice-punctured spheres, which will justify the notation.

By the first step we can assume that
$2X$ appears in a neighbourhood of $2M_0$, which is
controlled by Dehn surgery coefficients: For $n \in \N$ large enough we can
use Thurston's hyperbolic Dehn surgery theorem for $2M_0$ to produce a complete
hyperbolic orbifold of finite volume with the same topological type as $2X$ and cone-angles
$1/n$ around all edges. We can then continuously decrease the cone-angles of
$2X$ to the orbifold angles using the first step and obtain via Mostow-Prasad
rigidity for complete hyperbolic orbifolds of finite volume, cf.~\cite{Mos}, \cite{Pra}, that we have reached the above hyperbolic orbifold
structure at the end of the deformation.

We again use Thurston's hyperbolic Dehn surgery theorem for $2M_0$ to connect
$2M_0$ with $2X$ by a continuous family of hyperbolic
cone-3-manifolds $(D_{\tau})_{\tau \in (0,1]}$ with cone-angles
$(\tau\alpha_1, \ldots, \tau\alpha_N)$ and singular locus a link. Since we in
particular have local rigidity for hyperbolic cone-3-manifolds with totally
geodesic turnover boundary, we conclude that $D_{\tau}=2X_{\tau}$ for $\tau \in(0,1]$ with $X_1 =X$. It follows that for $\tau \rightarrow 0$ the
$X_{\tau}$ converge geometrically to a complete hyperbolic 3-manifold of finite volume $M_0$ with totally geodesic boundary consisting of thrice punctured spheres.

Mostow-Prasad rigidity, cf.~\cite{Mos}, \cite{Pra}, implies that the complete
hyperbolic structure on $2M_0$ is unique. By doubling we obtain that the
complete hyperbolic structure with totally geodesic boundary on $M_0$ is
unique. Furthermore we see that $2M_0$ is indeed the double of $M_0$ along its
totally geodesic boundary. In case the singular locus doesn't contain
vertices, essentially the same argument without doubling applies, cf.~\cite{Koj}.\\
\\
This finishes the proof of Theorem \ref{global_rigidity_hyperbolic} and
Theorem \ref{global_rigidity_hyperbolic_finite_volume}: 

Given a hyperbolic
cone-3-manifold $X$ as in Theorem \ref{global_rigidity_hyperbolic} or as in Theorem \ref{global_rigidity_hyperbolic_finite_volume},
we can construct by using Step 1 and Step 2 a family of hyperbolic
cone-3-manifolds $(X_{\tau})_{\tau \in (0,1]}$ of constant topological type
  and prescribed cone-angles connecting $X=X_1$ with a complete non-singular
  hyperbolic 3-manifold $M_0$ of finite volume, possibly with totally geodesic
  boundary consisting of thrice-punctured spheres. Since the complete hyperbolic
  structure with totally geodesic boundary on $M_0$ is unique,
  we conclude by local rigidity (resp.~a version of Thurston's hyperbolic
  Dehn surgery theorem applying to $M_0$, cf.~Appendix B in \cite{BP}) that the whole family $X_{\tau}$ is
  unique, in particular the initial hyperbolic cone-manifold structure on $X=X_1$.

\subsection{The spherical case}
Let $X$ be a compact orientable spherical cone-3-manifold with cone-angles
$\leq\pi$, which is not Seifert fibered. Let $N$ denote the number of edges contained in $\Sigma$ and let $(\alpha_1,
\ldots, \alpha_N)$ be the vector of cone-angles of $X$. We consider the interval
\begin{align*}
\mathcal{I}= \{ & t \in [0,1] : \exists \text{ a family of spherical cone-3-manifolds } (X_{\tau})_{\tau \in [0,t]}\\
& \text{of constant topological type with } X_0=X \text{ and cone-angles}\\
& ((1-\tau)\alpha_1+\tau\pi, \ldots, (1-\tau)\alpha_N+\tau\pi) \} \,.
\end{align*}
$\mathcal{I}$ is nonempty since clearly $0 \in
\mathcal{I}$. We claim that $\mathcal{I}$ is open and closed.\\
\\
{\bf(openness)} Let $t \in \mathcal{I}$ and a corresponding family of
spherical cone-3-manifolds $(X_{\tau})_{\tau \in [0,t]}$ be given. If $t<1$, we can further continuously increase the cone-angles of $X_t$ using local rigidity for compact spherical cone-3-manifolds, cf.~\cite{Wei}. For this result to hold we need the additional hypothesis that $X$, resp.~$X_t$ is not Seifert fibered.\hfill q.e.d.\\
\\
{\bf(closedness)} Suppose we are given a continuous family
$(X_{\tau})_{\tau\in[0,t_{\infty})}$ of spherical
  cone-3-manifolds. Consider a sequence $t_n \nearrow t_{\infty}$ and set
  $X_n:=X_{t_n}$. Schl\"afli's formula in the spherical case implies that
  $vol(X_{\tau})$ is a (strictly) increasing function of $\tau$, i.e.~$vol(X_{\tau}) \geq
  vol(X_0) > 0$. In the spherical case this is enough to conclude thickness,
  cf.~Section 9.3 in \cite{BLP},
  i.e.~there exists a $\rho$-thick basepoint $x_n \in X_n$ for some
  $\rho>0$. We may now invoke Theorem \ref{compactness} to extract a
  geometrically convergent subsequence. 

On the other hand, since the $X_n$ are Alexandrov spaces
  with curvature bounded below by $1$, we have $diam(X_n) \leq \pi$ and
  therefore a
  geometric limit will be compact. This implies that we can continuously
  extend the family $(X_{\tau})_{\tau \in [0,t_{\infty})}$ to the closed interval $[0,t_{\infty}]$.\hfill q.e.d.\\ 
\\
We have therefore established the existence of a continuous, angle-increasing
family of spherical cone-3-manifolds connecting $X=X_0$ with a spherical
orbifold $X_1$ having cone-angle $\pi$ around each edge. Global rigidity now
  follows from the following result of G.~de Rham:
\begin{thm}{\rm \cite{deR}}
A spherical structure on a closed, orientable 3-orbifold is
unique up to isometry.
\end{thm}
For a proof see \cite{deR} and \cite{Rot}. This finishes the proof of Theorem
\ref{global_rigidity_spherical}.


\end{document}